\documentclass[12pt]{article}
\usepackage{amsfonts}
\usepackage{amsmath}
\usepackage{amssymb}
\usepackage{color}
\font\smallit=cmti10  

\usepackage{bm}
\usepackage{amssymb,latexsym}
\usepackage{latexsym, amssymb, amsmath, amscd, amsfonts, graphicx} 
\usepackage{pifont}
\usepackage{CJK, psfrag,eepic}
\usepackage{mathrsfs}
\usepackage{txfonts}
\usepackage[perpage,symbol*]{footmisc}
\makeatletter

\renewcommand\section{\@startsection {section}{1}{\z@}
{-30pt \@plus -1ex \@minus -.2ex} {2.3ex \@plus.2ex}
{\normalfont\normalsize\bfseries}}

\renewcommand\subsection{\@startsection{subsection}{2}{\z@}
{-3.25ex\@plus -1ex \@minus -.2ex} {1.5ex \@plus .2ex}
{\normalfont\normalsize\bfseries}}

\renewcommand{\@seccntformat}[1]{\csname the#1\endcsname. }

\makeatother

\newtheorem{thm}{Theorem}[section]

\newtheorem{conj}[thm]{Conjecture}
\newtheorem{lem}[thm]{Lemma}

\newcommand{\N}{\mathbb N}
\newcommand{\Z}{\mathbb Z}
\newcommand{\R}{\mathbb R}

 \DeclareMathOperator{\ord}{ord}

 \DeclareMathOperator{\supp}{supp}

\def\pf{\noindent{\it Proof.} }
\def\qed{\nopagebreak \hfill $\Box$\medbreak}
\renewcommand{\b}{\hm\cdot}

\begin{document}

\title{\bf An Upper Bound for Davenport Constant of Finite Groups}

\maketitle

\begin{center}
 \vskip 20pt
{\bf  Weidong Gao$^1$,\ Yuanlin Li$^2$,\ Jiangtao Peng$^3$
\\ {\smallit $^1$Center for Combinatorics, LPMC-TJKLC,
Nankai University, Tianjin 300071, P.R. China}\\
{\smallit $^2$Department of Mathematics, Brock University, St. Catharines, Ontario, Canada L2S 3A1}\\
{\smallit $^3$College of Science, Civil Aviation University of China, Tianjin 300300, P.R. China}}\\

\footnotetext{E-mail address: wdgao@nankai.edu.cn (W.D. Gao),
yli@brocku.ca(Y.L. Li),  jtpeng@cauc.edu.cn (J.T. Peng)}
\end{center}

\begin{abstract}
Let $G$ be a finite (not necessarily abelian) group and let $p=p(G)$ be the smallest prime number
dividing $|G|$. We prove that $\mathsf d(G)\leq
\frac{|G|}{p}+9p^2-10p,$ where $\mathsf d(G)$ denotes the small
Davenport constant of $G$ which is defined as the maximal integer
$\ell$ such that there is a sequence over $G$ of length $\ell$
contains no nonempty one-product subsequence.
\end{abstract}


\section{Introduction}
Let $G$ be a finite group written multiplicatively. By a sequence $S$
over G, we mean a finite sequence of terms from $G$ which is unordered
and repetition of terms is allowed. We say that $S$ is an one-product
sequence if its terms can be ordered so that their product equals 1,
the identity element of the group. An one-product sequence $S$ is called a minimal one-product sequence if it cannot be partitioned into two nonempty, one-product
subsequences. The small Davenport constant $\mathsf d(G)$ is the maximal integer $t$ such that there is a sequence over
$G$ of length $t$ which contains no nonempty one-product
subsequence. The large Davenport constant $\mathsf D(G$) is the maximal
length of all minimal one-product sequences.  A simple argument \cite[Lemma 2.4]{GG2013} shows that
\begin{equation} \label{eq0.1}
\mathsf d(G) + 1 \leq \mathsf D(G)\leq |G|. \end{equation}
with equality in the
first bound when G is abelian, and equality in the second when $G$
is cyclic. The study of $\mathsf D(G) (= \mathsf d(G)+1)$, for $G$ abelian, is a
classical and very difficult problem in Combinatorial Number Theory.
When $G$ is non-abelian, there is more than one way to naturally
extend the definition of the Davenport constant. This was first done
by Olson and White \cite{OW} who introduced the small Davenport
constant $\mathsf d(G)$ and gave the general upper bound $\mathsf d(G)\leq
\frac{1}{2} |G|$ (for G non-cyclic) that was observed to be tight
for non-cyclic groups having a cyclic, index 2 subgroup.  When $G$ is a $p$-group,  $\mathsf d(G)$
was studied in \cite[Lemma 1.4]{Bass2007} and \cite{Dim2004}.
The large Davenport constant was introduced recently and studied in
\cite{GG2013} and \cite{Gryn2013}. A most recent result of
Grynkiewicz \cite{Gryn2013} states that $\mathsf d(G)+1\leq \mathsf D(G)\leq
\frac{2|G|}{p}$.  For an arbitrary finite non-abelian group $G$, let $p=p(G)$ denote the
smallest prime divisor of $|G|$. In this paper we  provide a better upper bound for the small Davenport constant and our main result is as follows.

\begin{thm}\label{theorem1}
Let $G$ be a finite noncyclic group of order $n$ and let $p$ be the
smallest prime divisor of $n$. Then $\mathsf d(G) \le \frac{n}{p}+
9p^2 -10p. $
\end{thm}

If $G$ has a cyclic subgroup $H$ of order $\frac{n}{p}$, then $H$ is
a normal subgroup of $G$ (\cite[Theorem 1]{Lam}).  Let $h$ be a generator of $H$ and
let $g\in G\setminus H$. Then the sequence
$$
S=\underbrace{g  \bm \cdot g\hm \cdot \ldots \hm \cdot
g}_{p-1} \hm \cdot \underbrace {h\hm \cdot  \ldots \hm \cdot
h}_{\frac{n}{p}-1}
$$
is an one-product free sequence of length $|S|=\frac{n}{p}+p-2.$
Therefore,
$$
\mathsf d(G)\geq \frac{n}{p}+p-2
$$
for any groups $G$ having a cyclic subgroup of order $\frac{n}{p}.$

We believe that the above mentioned lower bound is also an upper bound for the small Davenport constant.

\begin{conj}
Let $G$ be a finite noncyclic group of order $n$, and let $p$ be the
smallest prime divisor of $n$. Then $\mathsf d(G) \le \frac{n}{p}+
p-2. $
\end{conj}

\section{Preliminaries}

We use the notation and conventions described in detail in \cite{GG2013}.

For real numbers $a,b \in \R$, we set $[a,b] =\{x\in \Z: a \le x \le b\}$. If $A$ and $B$ are sets, we define the product-set as $AB=\{ab: a \in A, b\in B\}.$

Let $G$ be a finite multiplicative group. If $A \subseteq G$ is a nonempty subset, then denote by  $\langle A \rangle $  the subgroup  of $G$ generated by $A$. Recall that by a sequence over a group $G$ we mean a finite, unordered sequence where the repetition of elements is allowed. We view sequences over $G$ as elements of the free abelian monoid $\mathcal F (G)$ and we denote multiplication in $\mathcal F (G)$ by the bold symbol   $\bm \cdot$  rather than by juxtaposition and use brackets for all exponentiation in $\mathcal F (G)$.

A sequence $S \in \mathcal F(G)$ can be written in the form $S= g_1  \bm \cdot g_2 \bm \cdot \ldots \bm\cdot g_{\ell},$ where $|S|= \ell$ is the {\it length} of $S$. For $g \in G$, let
\begin{itemize}
\item[$\bullet$] $\mathsf v_g(S) = |\{ i\in [1, \ell] : g_i =g \}|\, \, $ denote the {\it multiplicity} of $g$ in $S$;
\item[$\bullet$] $\mathsf h(S) = \max \{ \mathsf v_g(S) : g \in G \}\, \, $ denote the {\it maximum multiplicity} of a term of  $S$;
\item[$\bullet$] $\supp(S) =\{g: \mathsf v_g(S) >0\}\,\, $ denote the {\it support} of $S$.
\end{itemize}
A sequence $T \in \mathcal F(G)$ is called a {\it subsequence } of $S$ and is denoted by $T \mid S$ if  $\mathsf v_g(T) \le \mathsf v_g(S)$ for all $g\in G$. Denote by   $T^ {[-1]} \bm\cdot S$ or $S \bm\cdot T^ {[-1]}$  the subsequence of $S$ obtained by removing the terms of $T$ from $S$.

If $S_1, S_2 \in \mathcal F(G)$, then $S_1 \bm\cdot S_2 \in \mathcal F(G)$ denotes the sequence satisfying that $\mathsf v_g(S_1 \bm\cdot S_2) = \mathsf v_g(S_1 ) + \mathsf v_g( S_2)$ for all $g \in G$. For convenience we  write
\begin{center}
 $g^{[k]} = \underbrace{g \bm\cdot \ldots \bm\cdot g}_{k} \in \mathcal F(G)\quad$ and $\quad T^{[k]} = \underbrace{T \bm\cdot \ldots \bm\cdot T}_{k} \in \mathcal F(G),$
\end{center}
for $g \in G, \, T \in \mathcal F(G)$ and $k \in \N_0$. Let $ T^{[-k]} = (T^{[k]})^{[-1]}$. If $ S_1$ and $  S_2$ are two subsequences of $S \in \mathcal F(G)$, then let $\gcd(S_1, S_2)$ denote the largest subsequence $T$ of $S$ such that $T \mid S_1$ and $T \mid S_2$.

Suppose $S= g_1 \b g_2 \bm \cdot \ldots \bm\cdot g_{\ell} \in \mathcal F(G)$, let $$\pi (S) = \{g_{\tau(1)}\ldots g_{\tau(\ell)}: \tau \mbox{ a permutation of } [1, \ell] \} \subseteq G$$  denote the {\it set of products} of $S$. Let $$\Pi(S) = \cup_{1 \le i \le \ell} \cup _{T\mid S,\ |T| = i}\pi(T)$$ denote the {\it set of all subsequence products} of $S$. The sequence $S$ is called
\begin{itemize}
\item[$\bullet$] {\it squarefree} if $\mathsf v_g(S) \le 1$  for all $g \in G$;
\item[$\bullet$]  {\it one-product} if $1 \in \pi(S)$;
\item[$\bullet$] {\it one-product free} if $1\not\in \Pi(S)$;
\item[$\bullet$] {\it minimal one-product } if $1\in \pi(S)$ and $S$ cannot be factored into two nontrivial, one-product subsequences.
\end{itemize}

We call $$\mathsf D(G)= \sup\{\, |S|: S \in \mathcal F(G) \mbox{ is
minimal one-product }\} \in \N_0 \cup \{ \infty\}$$ the {\it Large
Davenport constant} of $G$, and
$$\mathsf d(G)= \sup\{\, |S|: S \in \mathcal F(G) \mbox{ is one-product free }\} \in \N_0 \cup \{ \infty\}$$ the {\it small Davenport constant} of $G$.

\bigskip

\begin{lem} \cite{Olson1975} \label{k_set}
Let $G$ be  a group and let $S$ be an one-product free
sequence over $G$ of length $k$.  Then $|\Pi(S)| \ge
\frac{1}{9}k^2$.
\end{lem}

\begin{lem}\cite{OW}\label{kemperman}
Suppose $A$ and $B$ are finite subsets of a group and $1 \in A \cap
B$. If $1= ab \, (a\in A, b  \in B)$ has no solution except $a=b=1$,
then $|AB| \ge |A|+|B|-1.$
\end{lem}

The proof of lemma \ref{kemperman} may be found in Kemperman
\cite{Kemp} and Lemma \ref{kemperman} implies the following lemma.

\begin{lem}\label{partition}
Let $G$ be  a group and let $S$ be an one-product free
sequence over $G$. If  $S=S_1 \bm\cdot S_2 \bm\cdot \ldots
\bm\cdot  S_t$, then $|\Pi(S)| \ge \sum_{i=1}^{t}(|\Pi(S_i)|)$.
\end{lem}

\begin{lem}\label{length}
Let $S$ be an  one-product free sequence over a group $G$. Then
$|\Pi(S)| \ge |S|$.
\end{lem}

Let $N$ be a subgroup of a finite group $G$. For any element $a \in
G$, let $\overline{a} = aN$. For any subset $A$ of $G$, let
$\overline{A}=\{\overline{a}: a \in A\}$. Clearly $|\overline{A}|
\le |A|$, and  the equality holds if and only if  no two
elements of $A$ are in the same left coset of $N$.

\begin{lem} \label{subproduct}
Let $N$ be a subgroup of a finite group $G$. Let $A$ and $B$ be two
nonempty subsets of $G$ with $\overline{1} \in \overline{A} \cap
\overline{B}$. If $|\overline{B}| \geq 2$, then $|\overline{AB \cup BA}| \ge
\min \{p(G), |\overline{A}|+1\}$, where $AB=\{ab|a\in A, b\in B\}$, $BA=\{ba|a\in A, b\in B\}$
and $p(G)$ denotes the smallest prime divisor of $|G|$.
\end{lem}
\pf
Assume to the contrary that
\begin{align*}
|\overline{AB \cup BA}| \le \min \{p(G), \overline{A}+1\}-1 = \min
\{p(G)-1, |\overline{A}|\}.
\end{align*}
Then $|\overline{AB \cup BA}| \le |\overline{A}|$. Since $\overline{1} \in \overline{B}$, we infer that $\overline{A}
\subseteq \overline{AB} \subseteq \overline{AB \cup BA}$. Thus  $\overline{A} = \overline{AB \cup BA}$, so  $\overline{BA} \subseteq \overline{A} $. Since $\overline{1} \in \overline{A} \cap
\overline{B}$, we conclude that  $|\overline{BA} | \geq |\overline{A}| $. Thus we have $\overline{A}=\overline{BA}$.
Let  $b \in B \setminus N$ and let $a\in A$ such that $\overline{a} =\overline{1}$. Since $ \overline{BA} = \overline{A}$,  we obtain that
$\overline{b} = \overline{ba} \in \overline{BA}=\overline{A}$. Thus we have $\overline{b^2}=b\overline{b} \in \overline{BA} = \overline{A}$. Continuing this way, we obtain that $\overline{b^i} \in \overline{A}$  for all nonnegative integers $i$. Let $\ell$ be the smallest positive integer such $b^{\ell} \in N$. Then by the minimality of $\ell$, we get  $\ell \mid \ord(b)$. Since $b \not\in N$, we have $\ell >1$. Hence $\ell \ge p(G)$. Again, by the minimality of $\ell$, we conclude that $ \overline{1},  \overline{b}, \cdots, \overline{b^{\ell -1}}$ are distinct elements in $\overline{A}=\overline{AB \cup BA}$. So, $|\overline{AB \cup BA}| \ge p(G)$, yielding a
contradiction. \qed

\begin{lem}\label{coset}
Let $N$ be a subgroup of a finite group $G$, and let $S$ be a
sequence over $G\setminus N$. Then $|\{\overline{1}\} \cup
\overline{\Pi(S)}| \ge \min\{p(G), |S|+1\}$.
\end{lem}

\pf We proceed by induction on $|S|$.  If $|S|=1$ then
$|\{\overline{1}\} \cup \overline{\Pi(S)}|=2=|S|+1$. Assume that the
lemma is true for $|S|=k$ ($k\geq 1.$) and we want to prove it is
also true  for $|S|=k+1.$ Take any term $b|S$. Let $T=Sb^{-1}$. Then
$|T|=k$ and $|\{\overline{1}\} \cup \overline{\Pi(T)}| \ge
\min\{p(G), |T|+1\}$ by the inductive hypothesis. Let
$$
A=\{1\} \cup {\Pi(T)}
$$
and
$$
B=\{1, b\}.
$$
Then
$$
\overline{AB}\cup \overline{BA}\subseteq \{\overline{1}\} \cup
\overline{\Pi(S)}.
$$
It follows from Lemma \ref{subproduct} that $|\{\overline{1}\} \cup
\overline{\Pi(S)}|\geq |\overline{AB}\cup \overline{BA}|\geq \min
\{p(G), |A|+1\}\geq \min \{p(G), |T|+2\}=\min \{p(G), |S|+1\}.$ \qed

\section{Proof of the Main Theorem}

We are now ready to prove our main theorem.

{\bf Proof of Theorem \ref{theorem1}}.

Let $S \in \mathcal F(G)$ be
a sequence of length $\frac{n}{p}+c$ with $c=9p^2-10p+1$. Then we need to
show that $1 \in \Pi(S)$. Without loss of generality we may assume that
$\langle S \rangle =G$. We prove by the way of contradiction. Assume to the contrary that $S$ is
one-product free. Then
$$
|\Pi(S)|\le n-1.
$$

\noindent Let $t \in \N_0$ be maximal such that $S$ has a representation in
the form $S=S' \bm\cdot S_1 \bm\cdot S_2 \bm \cdot \ldots \bm \cdot
S_t$, where  $S_1, \ldots, S_t$ are squarefree, one-product free
subsequences of length $|S_\nu |=9p$ for all $\nu \in [1,t]$. Let
$d= |\supp (S')|$. Then
$$
|S'|+9pt=|S|=\frac{n}{p}+c.
$$
By the maximality of $t$ we get $0 \le d \le 9p-1.$
Since $S$ is one-product free, by Lemmas \ref{k_set}, \ref{partition},
 and \ref{length},  we have
\begin{align*}
\begin{array}{ll} n-1 & \ge |\Pi(S)| \ge |\Pi(S')|+ \sum_{i=1}^{t}(|\Pi(S_i)|) \\ &\ge
|S'|+9p^2t=p(\frac{n}{p}+c-|S'|)+|S'| \\& =n-1+p(c-|S'|)+|S'|+1.
\end{array}
\end{align*}
It follows that
\begin{equation} \label{eq1.1}
|S'|\geq c+1=9p^2-10p+2.
\end{equation}
Since $d\leq 9p-1$, we have
$$
\mathsf v_g(S')\geq p.
$$
for some $g\in G.$

For each $g \in G$ and each subsequence $T$ of $S$, let
$T_{\langle g \rangle}$ denote the subsequence of $T$ consisting of all terms in $\langle g \rangle$.
We first prove a useful claim.

{\bf Claim 1:} For each $g \in G$, let $C \mid S_{\langle g\rangle}$ and $D \mid S \b S_{\langle g\rangle}^{[-1]}$ with $|D| = p-1$. Then $|\Pi(C \b D)| \ge p |C|$.

 Let  $N = \langle
g\rangle$. Then $\Pi(C) \subseteq N$. By Lemma \ref{length} we have $|\Pi(C)| \ge |C|$. For any element $a \in G$, let $\overline{a} = aN$.
For any subset $A$ of $G$, let $\overline{A}=\{\overline{a}: a \in
A\}$. By Lemma \ref{coset} we have $|\{ \overline{1}\} \cup
\overline{\Pi(D)}| \ge p$.  Thus we obtain
that $|\Pi(C \b D)|\ge p |C|$. This proves our claim.

 We next rewrite $S'$ in a suitable form. Let $T=S'$ and choose $g_1 \in \supp(T)$. If $|T \b T^{[-1]}_{\langle g_1\rangle}| \geq p-1$, then $S'$ has a representation
$$S'=D_1 \b T_1  \b T',$$
where $T_1=T_{\langle g_1\rangle}$ and $D_1 | T \b T^{[-1]}_{\langle
g_1\rangle}$ with $|D_1|=p-1$. Let $T=T'$ and repeat the above
process on $T$. Thus $S'=D_1 \b T_1 \b D_2 \b T_2 \b T''$.
Continuing this way, we obtain that $S'$ has a representation
$$S'=D_1 \b T_1 \b D_2 \b T_2 \b \ldots \b
D_{\lambda} \b T_{\lambda} \b T,$$ where $\lambda=0$, or $\lambda
\in [1,d-1]$ and  $T_i | S'_{\langle g_i\rangle}, D_i | S' \b
S'^{[-1]}_{\langle g_i\rangle}$ with $|D_i|=p-1$, $\mathsf
v_{g_i}(T)=0$ for each $i\in [1,\lambda]$,  and $|T \b
T^{[-1]}_{\langle g\rangle}| \leq p-2$ for every $g \in \supp(T)$.

We now have the following two cases:

{\bf Case 1.} $|T_{\langle g\rangle}|\leq p-1$  for every $g \in
\supp(T)$. We note that $|T| \leq \min\{(p-1)(d-\lambda), 2p-3\}$.
Since $|S'|\geq 9p^2-10p+2 > 2p-3 $, we have $\lambda \geq 1$.
Therefore
$$
|T|\leq (p-1)(d-\lambda)\leq (p-1)(9p-2).
$$
 Note that
$$
\frac{n}{p}+c=|S|=9pt+(p-1)\lambda+|T_1|+\cdots +|T_{\lambda}|+|T|.
$$
By Lemmas \ref{partition}, \ref{length} and  Claim 1,  we conclude
\begin{align*}
\begin{array}{ll} n-1 & \ge |\Pi(S)| \ge |\Pi(S')|+ \sum_{i=1}^{t}(|\Pi(S_i)|) \\ &\ge
 |\Pi(T)| + \sum_{i=1}^{\lambda}  |\Pi(D_i\b T_i)| + \sum_{i=1}^{t}(|\Pi(S_i)|) \\ & \ge
|T|+p(|T_1|+\cdots +|T_{\lambda}|)+9p^2t \\
& =|T|+p(\frac{n}{p}+c-9pt-(p-1)\lambda-|T|)+9p^2t
\\& =n+p(c-(p-1)\lambda-|T|)+|T|  \ \ (\mbox{Since }0 \le |T| \le (p-1)(9p-2) )\\
&\geq n+p(c-(p-1)d)>n-1,
\end{array}
\end{align*}
yielding a contradiction.

{\bf Case 2.} There exists some element $g_{\lambda +1} \in \supp(T)$
 such that $|T_{\langle g_{\lambda+1} \rangle}|\geq p$.

If $T_1 \b \ldots \b T_{\lambda} \b T$ contains at least $p-1$ terms
not in $\langle g_{\lambda+1} \rangle$, then $S'$ has a
representation
$$S'=D_1 \b T_1' \b \ldots \b D_{\lambda} \b T_{\lambda}' \b D_{\lambda+1} \b T_{\lambda+1}$$
such that $T_i'\mid T_i$ for every $i\in [1,\lambda]$,
$T_{\lambda+1}=T_{\langle g_{\lambda+1}\rangle}$ and $D_{\lambda+1}$
is a sequence over $G\setminus \langle g_{\lambda+1}\rangle$ of
length $p-1$. As in Case 1, we get
\begin{align*}
\begin{array}{ll} n-1 & \ge |\Pi(S)| \ge |\Pi(D_{\lambda+1}\b T_{\lambda+1})+|\sum_{i=1}^{\lambda}  |\Pi(D_i\b T_i')| + \sum_{i=1}^{t}(|\Pi(S_i)|) \\ & \ge
p(|T_1'|+\cdots +|T_{\lambda}'|+|T_{\lambda+1}|)+9p^2t \\
& = p(\frac{n}{p}+c-9pt-(p-1)(\lambda+1))+9p^2t\\
&\geq n+p(c-(p-1)d)>n-1,
\end{array}
\end{align*}
yielding a contradiction. Therefore, we may assume that $T_1 \b
\ldots \b T_{\lambda} \b T$ contains at most $p-2$ terms not in
$\langle g_{\lambda+1} \rangle$. It follows that
$$
|S'_{\langle g_{\lambda+1}\rangle}|\geq |S'|-(p-2)-(|D_1|+\cdots
+|D_{\lambda}|)\geq |S'|-(p-2)-(9p-2)(p-1).
$$
If $S'$ contains at least $p-1$ terms not in $\langle g_{\lambda+1}
\rangle$, then $S'$ has a representation $S'=D_{\lambda+1} \b
T_{\langle g_{\lambda+1}\rangle} \b T'$ such that $D_{\lambda+1}$ is
a sequence over $G\setminus \langle g_{\lambda+1}\rangle$ of length
$p-1$ and $|T'|\leq (9p-3)(p-1)+p-2$. Now as in Case 1, we can
derive a contradiction.

Next we may assume that $S'$ contains at
most $p-2$ terms not in $\langle g_{\lambda+1} \rangle$. Let
$g=g_{\lambda+1}$. Then $S'$ has a representation
$$
S'=D \b S'_{\langle g \rangle}
$$
with $|D|\leq p-2.$ Note that $\langle S \rangle = G$ is noncyclic. Thus we infer that
$\langle g \rangle \ne G$. Since $S$ is one-product free,
$S_{\langle g \rangle}$ is one-product free, so we have $|S_{\langle
g \rangle}| \le |\langle g \rangle|-1 \le \frac{n}{p}-1$ and thus
$|S| - |S_{\langle g\rangle}|  \ge c > p$. Therefore,  $S$ contains
at least $ p-1$ terms not in $\langle g \rangle$. By renumbering if
necessary, we may assume that $S$ has a representation
$$
S=S_1 \b \ldots \b  S_{t'}\b D' \b S'_{\langle g \rangle} \b W
$$
with $t'\in [\max\{0,t-p+1+|D|\}, t]$, $D \mid D'$, ($D'$ is a sequence
over $G\setminus \langle g\rangle$ of length $p-1$) and $|W|\leq
(9p-1)(t-t')\leq (9p-1)(p-1-|D|)\leq (9p-1)(p-1)$. Again as in  Case 1, we can derive a contradiction.

In both cases we have found contradictions. Thus we must have $1 \in \Pi(S)$ and this completes the proof.

\qed

\bigskip
\noindent {\bf Acknowledgments}.

The research was carried out during a visit by the second author to
the Center for Combinatorics at Nankai University. He would like to gratefully acknowledge the kind hospitality from the host
institution. This work was supported in part by the 973
Program of China (Grant No. 2013CB834204),  the PCSIRT Project
of the Ministry of Science and Technology,  the National Science
Foundation of China and a Discovery Grant from the Natural Science and Engineering Research Council of Canada.

\end{document}